\documentclass{lmcs}
\pdfoutput=1

\usepackage{lastpage}
\lmcsdoi{17}{2}{3}
\lmcsheading{}{\pageref{LastPage}}{}{}%
{Nov.~01,~2019}{Apr.~14,~2021}{}

\usepackage[utf8]{inputenc}

\usepackage[numbers]{natbib}
\usepackage{amssymb,amsmath,amsthm}
\usepackage[inline]{enumitem}
\usepackage[matrix,arrow]{xy}

\newcommand{\iin}{\in}%
\newcommand{\LIN}{\mathsf{L}}
\newcommand{\A}{\mathfrak{A}}
\newcommand{\B}{\mathfrak{B}}
\newcommand{\D}{\mathfrak{D}}
\newcommand{\real}{\mathbb{R}}

\newcommand{\INT}[1]{\overline{#1}}
\newcommand{\rmd}{\mathrm{d}}

\newcommand{\sigmagen}{\mathbin{\widehat{\otimes}}}
\newcommand{\algen}{\otimes}

\newcommand{\KERNEL}{kernel}

\newcommand{\union}{\mathop{\textstyle\bigcup}}
\newcommand{\sm}{\setminus}
\newcommand{\spq}{\supseteq}
\newcommand{\sbq}{\subseteq}

\renewcommand{\phi}{\varphi}

\newcommand{\ep}{\varepsilon}

\newcommand{\N}{\mathbb{N}}

\newcommand{\defin}{\mathrel{\mathop:}=}

\newcommand{\calK}{\mathcal{K}}
\newcommand{\calD}{\mathcal{D}}

\begin{document}
\title{Semipullbacks of labelled Markov processes}
\author[Jan Pachl]{Jan Pachl\rsuper{a}}
\address{\lsuper{a}Department of Mathematics and Statistics, 
  York University,
  Toronto, Ontario, Canada.}
\author[Pedro S\'anchez Terraf]{Pedro S\'anchez Terraf\rsuper{b}}
\address{
  \lsuper{b}Universidad Nacional de C\'ordoba. 
  Facultad de Matem\'atica, Astronom\'{\i}a,  F\'{\i}sica y
  Computaci\'on. 
  Centro de Investigaci\'on y Estudios de Matem\'atica (CIEM-FaMAF),
  Conicet. C\'ordoba. Argentina.}
\thanks{Supported by Secyt-UNC project 33620180100465CB}
%
%% For Conicet:
%%
%% SEMIPULLBACK
%% RADON MEASURE
%% LABELLED MARKOV PROCESS
%% BISIMULATION
%% UNIVERSALLY MEASURABLE SET

\begin{abstract}
  A  \emph{labelled Markov process (LMP)} consists of a measurable
  space $S$ together with an indexed family of Markov kernels from $S$
  to itself. This structure has been used to model probabilistic
  computations in  Computer Science, and one of the main problems in
  the area is to define and decide whether two LMP $S$ and $S'$ ``behave
  the   same''. There are two natural categorical definitions of
  sameness of behavior: $S$  and $S'$  are \emph{bisimilar}
  if there exist  an LMP  $ T$ and measure preserving maps
  forming a diagram of the shape 
  $ S\leftarrow  T  \rightarrow{S'}$; and  they
  are \emph{behaviorally equivalent} 
  if there exist  some $ U$ and  maps forming a dual  diagram
   $ S\rightarrow  U  \leftarrow{S'}$.

  These two notions differ for general measurable spaces but Doberkat
  (extending  a result by Edalat)
  proved that they coincide for analytic Borel spaces, showing that
  from every 
  diagram  $ S\rightarrow  U  \leftarrow{S'}$ one can obtain a
  bisimilarity diagram as above. Moreover, the resulting square of
  measure preserving maps is commutative (a \emph{semipullback}).

  In this paper, we extend the previous result to measurable spaces $S$
  isomorphic to a  universally measurable subset of 
  a Polish space with the trace of the Borel  $\sigma$-algebra, using a
  version of Strassen's theorem on 
  common extensions of finitely additive measures.
\end{abstract}
\maketitle

\section{Introduction}
\label{sec:introduction}

Markov decision processes have been considered in the Computer
Science literature as a model for probabilistic computation. In this
context, a  \emph{labelled Markov process (LMP)} is
a structure $\mathbf{S} = (S,\Sigma, \{\tau_a  :  a\in L \})$ where $(S,\Sigma)$ is a
measurable space  and for $a\in L$, $\tau_a : S \times \Sigma
\rightarrow [0,1]$ is a \emph{Markov kernel}, i.e., a function  such
that for each fixed $s \in S$,  $\tau(s, \cdot)$ is a
finite positive measure bounded above by 1, and for each fixed $Q \in\Sigma$,
$\tau(\cdot, Q)$ is a $\Sigma$-$\B([0,1])$-measurable function. In one
interpretation of this computational model, the system
$\mathbf{S}$ stands at any particular time at a \emph{current state} $s_0\in S$, but
this information is hidden from the hypothetical \emph{users} of
$\mathbf{S}$, whose only interaction with the system is through
$L$. Intuitively, the user is presented
with a black box with buttons labelled by $L$, and a button $a$ is
available to be pressed whenever $\tau_a(s_0,S)>0$. A detailed
discussion of LMP and many motivating examples are to be found in
Desharnais' thesis~\cite{Desharnais}.

Of primary importance is to be able to determine when two such systems
$\mathbf{S}$ and $\mathbf{S'}$
behave the same way from the user viewpoint. That is, when a user
doing repeated experiments with $\mathbf{S}$ and $\mathbf{S'}$ would
conclude that they are indistinguishable.
Actually, for such
probabilistic systems there are at least two different ways to
formalize a notion of behavior, and they are intimately related to
measure-preserving maps.
\begin{defi}
  Let   $\mathbf S = (S,\Sigma, \{\tau_a : a\in L \})$ and  $\mathbf S' =
  (S',\Sigma', \{\tau_a' : a\in L \})$ be LMP. A \emph{zigzag
    morphism} $f:\mathbf{S}\to\mathbf{S'}$ is a surjective measurable map $f:
  (S,\Sigma)\to(S',\Sigma')$ such that for all $a\in L$ we have:
  \[\forall s\in S \;\forall Q\in\Sigma' : \tau_a(s,f^{-1}(Q)) =
  \tau_a'(f(s),Q).\]
\end{defi}
We say that $\mathbf S$  and $\mathbf{S'}$  are \emph{bisimilar}
if there exists  an LMP  $\mathbf T$ and zigzag morphisms
forming a diagram of the shape $\mathbf S\leftarrow \mathbf T
\rightarrow\mathbf{S'}$. This definition, in this categorical form, can be
traced to Joyal \emph{et al.} \cite{JOYAL1996164} and it provides one
of the possible formalizations of the concept of equality of
behavior. The second one is given by  the dual diagram: $\mathbf S$
and $\mathbf{S'}$  are \emph{behaviorally equivalent}
if there exists  an LMP  $\mathbf U$ and morphisms
forming a diagram of the shape $\mathbf S\rightarrow \mathbf U
\leftarrow\mathbf{S'}$. This notion, in turn, was introduced by Danos
\emph{et al.} \cite{coco}, and it can be
shown by functorial manipulations that behavioral equivalence (also
known as ``event bisimilarity'', i.e.\ the greatest \emph{event bisimulation}) is
a transitive relation in the category of LMP. It can be proved that
bisimilar LMP are behaviorally equivalent. Also, a neat
logical characterization of this last relation is given in \cite{coco}.
The originating papers of the concepts of LMP and bisimilarity with its 
logical characterization are 
Blute \emph{et al.} \cite{10.5555/788019.788852} and  Desharnais
\emph{et al.} \cite{10.5555/788020.788888}, and the presentation of
the results of both papers
were soon after streamlined in \cite{DEP}.
An alternative general source on the topic is Doberkat
\cite{doberkat2009stochastic}.

Some of the main problems in this area are to find conditions for the relation
of bisimilarity to be transitive, and more strongly, for behavioral
equivalence to entail bisimilarity. This is not true in the general
case \cite{Pedro20111048}, but there are various important positive
results which restrict or otherwise modify the category of processes and measurable spaces
considered.  The first one was obtained by Edalat \cite{Edalat} for a
category of LMP with a relaxed measurability condition on Markov kernels
(these are only required to be \emph{universally measurable}) over
analytic spaces: In such category of \emph{generalized} LMP, every
\emph{cospan}  $\mathbf S\rightarrow \mathbf U
\leftarrow\mathbf{S'}$ can be completed to a commutative square by
finding an appropriate $\mathbf{T}$ and arrows to $\mathbf{S}$,
$\mathbf{S'}$. This  $\mathbf{T}$ is called the \emph{semipullback} of
the cospan. Later, Doberkat
\cite{Doberkat:2005:SSR:1089905.1089907} obtained the same result now
properly for the category of LMP (with kernels as defined above) over
analytic spaces. He specifically showed the existence of
semipullbacks in the category of 
Markov kernels (that is, LMP with a singleton label set $L$) over
analytic state spaces and Borel zigzag maps; from
this, the result for general label sets follows.
%% It is also to be noted that the same result can be proved by combining
%% the logical characterization mentioned above with that obtained in
%% Blute \emph{et al.} \cite{10.5555/788019.788852}.

%
In the present paper we will show that the existence of semipullbacks
holds in the larger category of Markov kernels over universally
measurable spaces. Our proof does not rely on the existence of
disintegrations (regular conditional probabilities) as in \cite{Edalat}, but we use a result
about common extensions of finitely additive measures (Lemma~\ref{lem:extend}, a
version of Strassen's theorem). 
In Section~\ref{sec:statement-problem} we present
a related category, that of \emph{probability kernels}. The
main technical result of this paper is to show that this category has
semipullbacks. In Section~\ref{sec:preliminaries} we gather some
results on extensions of finitely additive
measures. Section~\ref{sec:construction} presents the construction of
the semipullback $S_3$ of a given cospan of probability kernels
$S_1\to S_0 \leftarrow S_2$; this is  essentially built over
the set-theoretic pullback of that diagram. The reduction of the
problem of Markov kernels and  general LMP to our result is done in
Section~\ref{sec:appl-probl-bisim};
in particular we show that LMP over coanalytic Borel spaces have
semipullbacks. We conclude with some counterexamples in the last section.

\section{Probability kernels}
\label{sec:statement-problem}
We find it technically convenient to describe the main construction
in terms of probability kernels from a fixed measurable space.
Let $(X,\Xi)$ and $(S,\Sigma)$ be two measurable spaces.
As in~\cite[Ch.1]{Kallenberg2002fmp},
a mapping $\mu\colon X\times \Sigma \to [0,\infty)$ is a
\emph{kernel from $X$ to $S$} if
$\mu(x,\cdot)$ is a measure on $\Sigma$ for each $x\in X$
and $\mu(\cdot, Q)$ is a $\Xi$--$\B([0,1])$-measurable function
on $X$ for each $Q\in\Sigma$.
From now on we write $\mu^x(Q)$ instead of $\mu(x,Q)$ for $x\in X$,
$Q\in \Sigma$.

Recall that a \emph{Radon measure} is a (non-negative) measure defined on the $\sigma$-algebra
of Borel sets of a topological space such that it is  inner regular with respect to compact sets; that
is, $\mu(B) = \sup \mu(K)$ where the supremum is over compact subsets
$K$ of $B$, for every Borel set $B$.
We say that $\mu$ is a \emph{probability kernel} if
$\mu^x(S)=1$ for all $x\in X$,
and a \emph{subprobability kernel} if
$\mu^x(S)\leq 1$ for all $x\in X$.
We say that $\mu$ is a \emph{Radon (sub)probability kernel}
if moreover $S$ is a topological space,
$\Sigma$ is its Borel $\sigma$-algebra
and every $\mu^x$ is a Radon (sub)probability measure on $\Sigma$.

Thus a Markov kernel in the definition of LMP above is a subprobability
kernel from $S$ to itself.

When $\mu$ is a kernel from $X$ to $S$, 
we write $(S,\Sigma,\mu)$ instead of $\mu$
when $(X,\Xi)$ is understood and we wish to make $S$ and $\Sigma$ explicit.

For a fixed $(X,\Xi)$, kernels from $X$ form a category
with surjective measure-preserving maps as morphisms:

\begin{defi}
  Let $(X,\Xi)$ be a fixed measurable space.
  For $j=1,2$ and $x\in X$ let $(S_j,\Sigma_j,\mu^x_j)$ be a measure space
  such that $\mu_j$ is a kernel from $X$ to $S_j$.
  A mapping $h\colon S_1\to S_2$ is a \emph{\KERNEL{} morphism from $\mu_1$ to $\mu_2$}
  if it is $\Sigma_1$--$\Sigma_2$ measurable, $h(S_1)=S_2$,
  and $\mu^x_1(h^{-1}(A)) = \mu^x_2(A)$ for all $x\in X$, $A\in\Sigma_2$.
  A morphism $h$ from $\mu_1$ to $\mu_2$ is sometimes written
  $h\colon (S_1,\Sigma_1,\mu_1)\to (S_2,\Sigma_2,\mu_2)$, 
  or simply $h\colon S_1\to S_2$ when $\Sigma_j$ and $\mu_j$ are understood.
\end{defi}
We find that notation is somewhat simpler when we work with kernel
morphisms rather than zigzag morphisms. Once we prove the existence of
a semipullback for kernel morphisms, the existence for zigzag
morphisms will easily follow (Section~\ref{sec:appl-probl-bisim}).

In the present paper we prove:
\begin{thm}
    \label{th:main}
  Let $(X,\Xi)$ be a fixed measurable space.
  Consider the category in which each object is 
  a Radon subprobability kernel from $X$ to a separable metric space,
  and morphisms are \KERNEL{} morphisms.
  Every cospan 
  $(S_1,\Sigma_1,\mu_1)\rightarrow (S_0,\Sigma_0,\mu_0) \leftarrow (S_2,\Sigma_2,\mu_2)$ 
  has a semipullback $(S_3,\Sigma_3,\mu_3)$ such that $S_3$ is the set pullback 
  of  $S_1 \rightarrow S_0 \leftarrow S_2$.
  Moreover $S_3$ is a measurable subset of $S_1\times S_2$.
\end{thm}

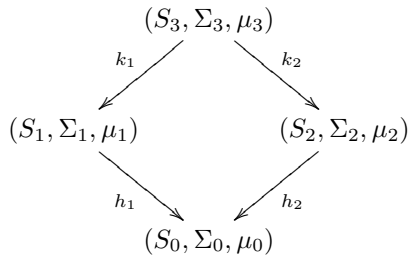
\begin{figure}[h]
  \[\xymatrix@C=-1ex@R=6ex{%
    & {(S_3,\Sigma_3,\mu_3)} \ar[dl]_{k_1} \ar[dr]^{k_2} \\
    {(S_1,\Sigma_1,\mu_1)} \ar[dr]_{h_1} & &
    {(S_2,\Sigma_2,\mu_2)} \ar[dl]^{h_2} \\
    &  {(S_0,\Sigma_0,\mu_0)}%
  }\]
  \caption{A semipullback.}
  \label{fig:diagram}
\end{figure}

We will first prove the theorem for probability kernels. For this we fix, up to Section~\ref{sec:construction}, three Radon probability
kernels $(S_j,\Sigma_j,\mu_j)$ (for $j=0,1,2$) as in the statement of
Theorem~\ref{th:main}, 
and for $j=1,2$, \KERNEL{} morphisms $h_j\colon S_j \to S_0$.  Our goal
is to construct a semipullback of $h_1$, $h_2$,
i.e. $(S_3,\Sigma_3,\mu_3)$ and for $j=1,2$, morphisms  $k_j\colon S_3 \to S_j$
such that $h_1\circ k_1 = h_2 \circ k_2$ (see Figure~\ref{fig:diagram}).

We first proceed to construct the pullback of the mappings $h_1$ and
$h_2$ in the category of measurable spaces, whose upper vertex will be the underlying
space of the semipullback.
Let $\pi_j\colon S_1\times S_2 \to S_j$, $j=1,2$, be the natural projections.
Denote by $\Sigma_1\sigmagen\Sigma_2$ the smallest $\sigma$-algebra on
$S_1\times S_2$ for which
$\pi_j$ are $\Sigma_j$-measurable.

Define
\begin{equation}\label{eq:1}
  \begin{split}
    S_3 & := \{ (x_1,x_2) \iin S_1\times S_2  :  h_1(x_1) = h_2(x_2) \} \\
    \Sigma_3 & := \{ A\cap S_3  :  A \iin \Sigma_1\sigmagen\Sigma_2 \} \\
    k_j & := \text{ restriction of } \pi_j \text{ to } S_3, \text{ for } j=1,2.
  \end{split}
\end{equation}
Then $k_j(S_3)=S_j$ for $j=1,2$ (because $h_j(S_j)=S_0$).

All that remains now is to construct the probability kernel $\mu_3$; this will be done in several steps.
We define a countable algebra $\A \subseteq\Sigma_3$ that generates $\Sigma_3$ as a $\sigma$-algebra,
and finitely additive measures $\nu^x_3$ on $\A$.
In defining $\nu^x_3$ we use a constructive variant of the Hahn-Banach theorem,
to ensure that $\nu^x_3$ is a measurable function of $x$.
Then we prove that $\nu^x_3$ is countably additive, so that it extends
to a countably additive measure $\mu^x_3$ on $\Sigma_3$.

\section{Preliminaries}
\label{sec:preliminaries}

In this section we establish notation and several results that will be needed in the main construction.

If $S$ is a set and $V$ is a set of real-valued  functions on $S$, 
write $V^+ := \{ f\iin V : f \geq 0 \}$.
Here $\geq$ is the pointwise partial order.

When $\B_1$ and $\B_2$ are algebras of subsets of $S_1$ and $S_2$,
we denote by $\B_1\algen\B_2$ the algebra of subsets of $S_1\times S_2$
consisting of finite unions of sets of the form $B_1\times B_2$, $B_j\in\B_j$ for $j=1,2$.
Recall that  we denote by $\B_1\sigmagen\B_2$ the $\sigma$-algebra
generated by $\B_1\algen\B_2$.

If $\A$ is an algebra of sets, denote by $\LIN(\A)$ the space of simple $\A$-measurable
functions;
that is, functions of the form $\sum_{i\in F} r_i \chi_{A_i}$ where
$F$ is a finite set,   $A_i\in\A$
and $r_i\iin\real$ for $i\iin F$.
If in addition $\nu$ is a finitely additive measure on $\A$, denote by $\INT{\nu}$ its integral
defined on $\LIN(\A)$:
\[
\INT{\nu} (f) := \int f \rmd \nu \quad \text{for } f \iin \LIN(\A).
\]
In other words,
$\INT{\nu} (f)= \sum_{i\in F} r_i \nu(A_i)$ when $f=\sum_{i\in F} r_i \chi_{A_i}$.
This is well-defined--- see e.g.~\cite[Ch.~III]{dunford-schwartz1}.

Note that if $\nu^x (A)$ is a measurable function of $x$ for every $A\iin\A$ then
so is $\INT{\nu^x} (f)$ for every $f\iin\LIN(\A)$.

\begin{lem}
  \label{lem:aimage}
  Any Borel-measurable image of a compact metrizable space in a
  separable metric space is analytic,
  and therefore universally measurable.
\end{lem}

\begin{proof}
  Let $K$ be compact, $S$ a separable metric space, and $h:K\to S$
  Borel. Then the  completion  $\hat S$ of $S$ is analytic and
  $h:K\to \hat S$ is also Borel. Hence
  by \cite[423G(b)]{fremlinMTvol4} the image is analytic, hence by
  \cite[434D(c)]{fremlinMTvol4}  it is universally measurable.
\end{proof}

\begin{lem}
  \label{lem:aquoient}
  Let $S$ and $S_0$ be separable metric spaces and $\Sigma$,
  $\Sigma_0$ their Borel $\sigma$-algebras.
  Let $\mu$ and $\mu_0$ be Radon probability measures on $\Sigma$ and
  $\Sigma_0$, respectively.
  Let $h\colon S \to S_0$ be a Borel-measurable measure-preserving mapping.
  Then for every $A\iin\Sigma$ there exists $D\iin\Sigma_0$ such that $D\subseteq h(A)$ and
  $\mu_0(D) \geq \mu(A)$.
\end{lem}

\begin{proof}
  Since $\mu$ is inner regular with respect to compact sets, we can
  find a countable family $\calK$ of compact subsets of $A$ such that
  \[
  \mu(A)=\sup\{\mu(K) : K\in\calK\}.
  \]
  For every  such
  $K$ the  image $h(K)$ is
  universally measurable by Lemma~\ref{lem:aimage}. Therefore, there
  exist $B_K, B'_K\in \Sigma_0$ such that $B_K\sbq h(K)\sbq B'_K$ and
  $\mu_0(B'_K\sm B_K) =0$. Since $h$ is measure-preserving, we have
  \[
  \mu_0(B_K)=\mu_0(B'_K)=\mu(h^{-1}(B'_K)).
  \]
  Now $h^{-1}[B'_K]\spq K$, and therefore
  $\mu_0(B_K)=\mu(h^{-1}(B'_K))\geq\mu(K)$ for each $K\in\calK$. Hence
  we may take $D\defin \union \{B_K :  K\in\calK\}$.
\end{proof}

\begin{lemC}[{\cite[457C]{fremlinMTvol4}}]
  \label{lem:extend}
  Let $\A$ be an algebra of subsets of a set $S$, and $\A_j$, $j=1,2$,
  two of its subalgebras.
  Let $\nu_j \colon \A_j \to [0,1]$, $j=1,2$, be finitely additive measures
  such that $\nu_1(S)=\nu_2(S)=1$.
  Assume that $\nu_1(A_1) + \nu_2(A_2) \leq 1$ whenever
  $A_j \iin \A_j$, $j=1,2,$ are such that $A_1 \cap A_2 = \emptyset$.
  Then there exists a finitely additive measure $\nu \colon \A \to [0,1]$
  that extends both $\nu_1$ and $\nu_2$. \qed
\end{lemC}

We need the following variant of the Hahn-Banach theorem,
which preserves measurability.

\begin{lem}
  \label{lem:measurable}
  Let $S$ be a non-empty set.
  Let $V$ be a linear space of bounded functions, and $W$ a subspace of
  $V$ that contains all constant functions on $S$. Assume also that V has a
  countable basis as a linear space.
  Let $\Psi^x \colon W \to \real$, $x\in X$, be a collection of linear functionals on $W$ such that
  $\Psi^x(f)$ is a measurable function of $x$ for every $f\iin W$, $\Psi^x(1)=1$ and
  $\Psi^x(f) \geq 0$ for $f\iin W^+$, $x\iin X$.
  Then there is a collection of linear functionals $\Phi^x \colon V \to \real$, $x\in X$,
  that extend $\Psi^x$
  and such that $\Phi^x(f)$ is a measurable function of $x$ for every $f\iin V$,
  and $\Phi^x(f) \geq 0$ for $f\iin V^+$, $x\iin X$.
\end{lem}

\begin{proof}
Extend $\Psi^x$ one dimension at a time.
Assume that $\Phi^x$ has been defined on a linear subspace $U\supseteq W$ and consider $f_0\iin V \setminus U$,
We are going to extend $\Phi^x$ to $U+\real f_0$ as follows.
\begin{align*}
p^x(f) & := \inf \{ \Phi^x(g) : g\iin U \text{ and } g \geq f \}  \quad\text{for } f \iin U+\real f_0\\
\Phi^x(f_0) & := \inf \{ p^x(g+f_0) - \Phi^x(g)  : g\iin U \}.
\end{align*}
Thus $p^x$ is subadditive and positively homogeneous on $U+\real f_0$.
We claim that $\Phi^x(f) \leq p^x(f)$ for every $f \iin U+\real f_0$.
To prove the claim, write $f=u+rf_0$ where $u\iin U$ and $r\iin\real$, and distinguish two cases: \\
For $r>0$ use $g=u/r$ in the definition of $\Phi^x(f_0)$ to get
\begin{align*}
\Phi^x(f) & = \Phi^x(u) + r\Phi^x(f_0) = \Phi^x(u) + \inf \{ r p^x(g+f_0) - r \Phi^x(g)  : g\iin U \}  \\
& \leq \Phi^x(u) + r p^x((u/r)+f_0) - r \Phi^x(u/r) = p^x(f).
\end{align*}
When $r<0$, for every $g\iin U$ we have
\[
\Phi^x(g) - \Phi^x(u/r)  = p^x(g - (u/r) )
\leq p^x(g + f_0 ) + p^x(- (u/r) - f_0)
\]
Therefore
\[
-p^x(- (u/r) - f_0) - \Phi^x(u/r)
 \leq \inf \{ p^x(g + f_0 ) - \Phi^x(g) : g\iin U \} = \Phi^x(f_0),
\]
and finally
\[
\Phi^x(f) = r\Phi^x(u/r) + r\Phi^x(f_0)
 \leq -r p^x(- (u/r) - f_0) = p^x(u + rf_0) = p^x(f).
\]
That proves the claim.
It follows that if $f\iin (U+\real f_0)^+$ then
\[
\Phi^x(f) = -\Phi^x(-f) \geq -p^x(-f) \geq 0.
\]

  To prove that $\Phi^x(f_0)$ is a measurable function of $x$,
  fix a countable basis $C$ of $U$ such that $1\iin C$ and define $\widetilde{U}$
  to be the set of finite linear combinations of elements of $C$ with rational coefficients.
  Then
  \begin{align*}
    p^x(f) & = \inf \{ \Phi^x(g) : g\iin \widetilde{U} \text{ and } g \geq f \}
    \quad\text{for } f \iin U+\real f_0\\
    \Phi^x(f_0) & = \inf \{ p^x(g+f_0) - \Phi^x(g)  : g\iin \widetilde{U} \}
  \end{align*}
  so that $x\mapsto \Phi^x(f_0)$ is the infimum of a countable set of measurable functions.
\end{proof}

The next lemma is a variant of a theorem of Marczewski and Ryll-Nardzewski~\cite{Marczewski1953}.
When $\B_j=\Sigma_j$, this is a special case of \cite[454C]{fremlinMTvol4}.

\begin{lem}
  \label{lem:compactproduct}
Let $S_1$ be a Hausdorff topological space, $\Sigma_1$ its Borel $\sigma$-algebra and
$\mu_1\colon\Sigma_1 \to[0,1]$ a Radon probability measure.
Let $(S_2,\Sigma_2,\mu_2)$ be any probability space.
Denote by $\pi_j\colon S_1 \times S_2 \to S_j$ the natural projections.
For $j=1,2$, let $\B_j \subseteq \Sigma_j$ be an algebra of subsets of $S_j$.
Let $\mu\colon\B_1\algen\B_2 \to [0,1]$ be a finitely additive measure such that
$\mu(\pi_j^{-1}(B_j))=\mu_j(B_j)$ for $j=1,2$ and all $B_j\in\B_j$.
Then $\mu$ has an extension
to a countably additive measure on the $\sigma$-algebra $\B_1\sigmagen\B_2$.
\end{lem}

\begin{proof}
This is a minor modification of the proof of \cite[454C]{fremlinMTvol4}.
Let $\D$ be the set of finite unions of sets of the form $C\times B_2$
where $C$ is a compact subset of $S_1$ and $B_2\in\B_2$.
As $\mu_1$ is a Radon measure,
it follows that for every $\ep>0$ and $B\in\B_1\algen\B_2$ there are $D\in\D$ and $E\in\Sigma_1$
such that $D\subseteq B$, $\mu_1(E)<\ep$ and $B\subseteq D\cup(E\times S_2)$.

Now let $\{B_i\}_{i\in\N}$ be a non-increasing sequence of sets
in $\B_1\algen\B_2$ with empty intersection.
To prove that $\lim_i \mu(B_i) = 0$, take any $\ep>0$.
There are $D_i^\prime\in\D$ and $E_i^\prime\in\Sigma_1$ such that
$D_i^\prime\subseteq B_i$, $\mu_1(E_i^\prime)<2^{-i}\ep$
and $B_i\subseteq D_i^\prime\cup(E_i^\prime\times S_2)$.
Set $D_n := \bigcap_{i\leq n} D_i^\prime$ and $E_n := \bigcup_{i\leq n} E_i^\prime$
for each $n$.
Then $\{D_n\}_n$ is a non-increasing sequence of sets in $\D$,
$D_n\subseteq B_n$, $\mu_1(E_n)<2\ep$ and
$B_n\subseteq D_n\cup(E_n\times S_2)$.

For $n\in\N$ and $y\in S_2$ set
$D_n^y := \pi_1( D_n \cap \pi_2^{-1}(y) ) $
and $H_n := \pi_2(D_n)$.
Then $\{D_n^y\}_n$ and $\{H_n\}_n$ are non-increasing sequences of subsets of $S_1$ and $S_2$,
respectively.
The sets $D_n^y$ are compact and $H_n\in\B_2$.
Next $\bigcap_n D_n^y = \emptyset$ because $\bigcap_n D_n \subseteq \bigcap_n B_n = \emptyset$.
Hence for every $y\in S_2$ there is $n$ such that $D_n^y = \emptyset$,
which means that $\bigcap_n H_n = \emptyset$.
It follows that
\begin{multline*}
  \lim_i \mu(B_i)
  \leq  \lim_i \mu(S_1\times H_i)  + \lim_i \mu( B_i \setminus
  (S_1\times H_i) )\\ 
  \leq \lim_i \mu_2(H_i) + \lim_i \mu_1(E_i) \leq 2 \ep.
\end{multline*}
We have proved that $\lim_i \mu(B_i) = 0$.
By~\cite[413K]{fremlinMTvol4}
$\mu$ has an extension
to a countably additive measure on the $\sigma$-algebra generated by $\B_1\algen\B_2$.
\end{proof}

\section{Proof of the main theorem}
\label{sec:construction}

In this section we complete the proof of Theorem~\ref{th:main}.
Recall that $S_j$, $j=0,1,2$, are separable metric spaces,
$\Sigma_j$ are their Borel $\sigma$-algebras,
and $\mu_j$ are Radon probability kernels from $X$ to $S_j$.

Our goal is to construct a semipullback $(S_3,\Sigma_3,\mu_3)$ of the
cospan   $(S_1,\Sigma_1,\mu_1)\rightarrow (S_0,\Sigma_0,\mu_0) \leftarrow (S_2,\Sigma_2,\mu_2)$.
In Section~\ref{sec:statement-problem} we have already defined $S_3$,
$\Sigma_3$, and the maps $k_1$ and
$k_2$ closing the diagram. To complete the construction, we will define the measures $\nu^x_3$
on $\A$ by using the results of the previous section.
\medskip

For $j=0,1,2$, fix countable algebras $\B_j\subseteq\Sigma_j$ such that
\begin{itemize}
\item
$\B_j$ generates $\Sigma_j$ as a $\sigma$-algebra for $j=0,1,2,$ and
\item
$h_j^{-1}(B_0)\in\B_j$ whenever $B_0\in\B_0$, for $j=1,2$.
\end{itemize}
For $j=1,2$, define $\A_j := \{ k_j^{-1} (B)  :  B\iin\B_j \}$.
Let $\A$ be the algebra of subsets of $S_3$ generated by $\A_1 \cup \A_2$.
Then $\A$ is countable and it generates $\Sigma_3$ as a $\sigma$-algebra.

For $j=1,2$, and $B\iin\B_j$, let $\nu^x_j (k_j^{-1} (B)):= \mu^x_j (B)$.
As $k_j(S_3)=S_j$, this is well defined
and $\nu^x_j$ is a finitely additive measure on $\A_j$.

Take any $A_j \iin \A_j$, $j=1,2,$ such that $A_1 \cap A_2 = \emptyset$.
Then $A_j=k_j^{-1} (B_j)$ for some $B_j\iin\B_j$.
By Lemma~\ref{lem:aquoient} there are $D_j\iin\Sigma_0$ such that
$D_j\subseteq h_j(B_j)$ and $\mu^x_0(D_j) \geq \mu^x_j(B_j)$.
From the definition of $S_3$ we get $h_1(B_1) \cap h_2(B_2) = \emptyset$,
hence $D_1 \cap D_2 = \emptyset$.
Therefore
\[
\nu^x_1(A_1) + \nu^x_2(A_2) = \mu^x_1(B_1) + \mu^x_2(B_2)
\leq \mu^x_0(D_1) + \mu^x_0(D_2) \leq 1.
\]
By Lemma~\ref{lem:extend} there is a finitely additive measure $\nu^x$
on $\A$ that extends both $\nu^x_1$ and $\nu^x_2$.
As the proof of Lemma~\ref{lem:extend} relies on the axiom of choice,
$\nu^x(A)$ is not necessarily a measurable function of $x$ for every $A\iin\A$.
However, observe that $\INT{\nu^x}(f)$ is a measurable function of $x$
for every $f\iin \LIN(\A_1) + \LIN(\A_2)$.
Indeed, if $f=f_1+f_2$, $f_j\iin\LIN(\A_j)$ for $j=1,2$,
then
\[
\INT{\nu^x}(f) = \INT{\nu^x_1}(f_1) + \INT{\nu^x_2}(f_2)
\]
by the linearity of integral, so that $\INT{\nu^x}(f)$ is a sum of two measurable functions of $x$.

By Lemma~\ref{lem:measurable} with $W=\LIN(\A_1) + \LIN(\A_2)$ there is a linear functional $\Phi^x \colon \LIN(\A) \to \real$
that agrees with $\INT{\nu^x}$ on $\LIN(\A_1) + \LIN(\A_2)$ and such that
$\Phi^x(f)$ is a measurable function of $x$ for every $f\iin\LIN(\A)$,
and $\Phi^x(f) \geq 0$ for $f\iin \LIN(\A)^+$, $x\iin X$.
Now $\nu^x_3(A):=\Phi^x(\chi_A)$ for $A\iin\A$ defines a finitely additive measure $\nu^x_3 \geq 0$
on $\A$ such that $\nu^x_3(S_3)=1$ and $x\mapsto \nu^x_3(A)$ is measurable for every $A\in\A$.

We have $\A = \{ B \cap S_3 : B \in \B_1\algen\B_2 \}$.
For $B\in\B_1\algen\B_2$ define $\mu^x (B) := \nu^x_3 (B \cap S_3 )$.
By Lemma~\ref{lem:compactproduct} each $\mu^x$ extends to a countably additive
measure $\widehat{\mu}^x$ on the $\sigma$-algebra
$\B_1\sigmagen\B_2 = \Sigma_1\sigmagen\Sigma_2$, which is the Borel $\sigma$-algebra of the product topology on $S_1\times S_2$.
By~\cite[454A(a)]{fremlinMTvol4}, $\widehat{\mu}^x$ is a Radon measure.

\begin{lem}
  \label{lem:definability-S_3}
  $S_3\in \Sigma_1\sigmagen\Sigma_2$; moreover, we have
  \begin{equation}\label{eq:definability-S_3}
    S_1\times S_2 \setminus S_3
    = \bigcup \{ h_1^{-1}(B_0) \times h_2^{-1}(S_0 \setminus B_0) : B_0 \in \B_0 \}
  \end{equation}
  and $\widehat{\mu}^x (S_3) = 1$.
\end{lem}
\begin{proof}
  Take any $(x_1,x_2) \in S_1\times S_2 \setminus S_3$.
  Then $h_1(x_1)\neq h_2(x_2)$,
  hence there is $B_0\in\B_0$ such that $h_1(x_1)\in B_0$
  and $h_2(x_2)\not\in B_0$,
  which means $(x_1,x_2)\in h_1^{-1}(B_0) \times h_2^{-1}(S_0 \setminus B_0)$.

  Each set $B := h_1^{-1}(B_0) \times h_2^{-1}(S_0 \setminus B_0)$,
  where $B_0\in\B_0$,
  is in the algebra $\B_1\algen\B_2$ and $\widehat{\mu}^x(B) =
  \mu^x(B) = \nu^x_3(\emptyset) = 0$.
  It follows that $S_3\in\Sigma_1\sigmagen\Sigma_2$ and
  $\widehat{\mu}^x (S_3) = 1$.
\end{proof}

By Lemma~\ref{lem:definability-S_3}, $S_3$ is a measurable subset of $S_1\times S_2$.
Define $\mu_3^x$ to be the restriction of $\widehat{\mu}^x$ to the $\sigma$-algebra $\Sigma_3$.

It remains to be proved that for every $E\in\Sigma_3$ the function $x\mapsto\mu_3^x$
is measurable.
To that end define
\[
  \calD \defin \{ E\in \Sigma_3 : x\mapsto\mu_3^x(E)\text{ is measurable}\}.
\]
Then $\A\subseteq\calD$, and $\calD$ is closed under complements and unions of disjoint sequences.
By the Monotone Class Theorem \cite[136B]{fremlinMTvol1} we have $\calD=\Sigma_3$.

That completes the proof of Theorem~\ref{th:main} for the case of
probability kernels. To extend the result to subprobability kernels we
work as follows.  Let $(S_1,\Sigma_1,\mu_1)\rightarrow (S_0,\Sigma_0,\mu_0) \leftarrow (S_2,\Sigma_2,\mu_2)$
be a cospan of Radon subprobability kernels, 
where $S_j$ are separable metric spaces. 

Define $\bar S_j \defin S_j \oplus \{s_j\}$ where $s_j\notin S_j$ for
each $j=0,1,2$ and for measurable $E\sbq \bar S_j$, let
\[
\bar\mu^x_j(E)\defin \mu_j^x(E \cap S_j) + (1-\mu_j^x(S_j))\cdot\chi_E(s_j) 
\]
Then $\bar\mu_j$ are Radon probability kernels.
We also extend the maps $h_j$ by stipulating
\[
\bar h_j(x)\defin
\begin{cases}
  h_j(x)& x\neq s_j\\
  s_0 & x=s_j
\end{cases}
\]
for $j=1,2$. 
Then $\bar h_j$ are \KERNEL{} morphisms.

By Theorem~\ref{th:main} for probability kernels, the cospan
$ \bar S_1 \rightarrow \bar S_0 \leftarrow \bar S_2$
has a semipullback $(\bar S_3,\bar \Sigma_3,\bar \mu_3)$ with \KERNEL{} morphisms
$k_j\colon\bar S_3 \to \bar S_j$,
and $\bar S_3 \subseteq \bar S_1 \times \bar S_2$ is the set pullback.
Hence $\bar S_3 = S_3\oplus\{(s_1,s_2)\}$ where
$S_3$ is the set pullback of $S_1 \rightarrow S_0 \leftarrow S_2$.

We can take $\mu^x_3$ to be the restriction of $\bar \mu^x_3$ 
to $\Sigma_3 \defin \{ E \cap S_3 : E \in \bar \Sigma_3 \}$. 
It is straightforward to
check that the restrictions $k_j\restriction S_3$ are \KERNEL{} morphisms
from $S_3$ onto $S_j$ for $j=1,2$, and we are done.

\section{Application to the problem of bisimulation}
\label{sec:appl-probl-bisim}

\subsection{Labelled Markov processes with Radon measures}
\label{sec:markov-kernels}

By Theorem~\ref{th:main}, semipullbacks exist in a certain category of subprobability kernels
from a fixed measurable space $(X,\Xi)$.

As a corollary we obtain the following theorem,
which asserts the existence of semipullbacks in the corresponding category
of LMP and zigzag morphisms:

\begin{thm}
    \label{th:Markov}
  Consider the category in which objects are LMP $(S,\Sigma, \{\tau_a : a\in L \})$
  such that $S$ is a separable metric space and $\tau_a(s,\cdot)$ are Radon measures,
  with zigzag morphisms.
  In this category every cospan has a semipullback.
  
  Moreover, every cospan
  \[
  (S_1,\Sigma_1,\{\tau_{1a} : a\in L \})
  \rightarrow (S_0,\Sigma_0,\{\tau_{0a} : a\in L \}) 
  \leftarrow (S_2,\Sigma_2,\{\tau_{2a} : a\in L \})
  \]
  has a semipullback $(S_3,\Sigma_3,\{\tau_{3a} : a\in L \})$ such that $S_3$ is the set pullback
  of  $S_1 \rightarrow S_0 \leftarrow S_2$
  and $S_3$ is a measurable subset of $S_1\times S_2$.
\end{thm}

\begin{proof}
First we deal with the LMP for which the label set $L$ has a single element $a$,
and write $\tau=\tau_a$.

Let $(S_1,\Sigma_1,\tau_1)\rightarrow (S_0,\Sigma_0,\tau_0) \leftarrow (S_2,\Sigma_2,\tau_2)$
be a cospan in the given category,
with connecting zigzags $h_j\colon S_j  \to S_0$, $j=1,2$.
As in Theorem~\ref{th:main}, take the measurable pullback $(S_3,\Sigma_3)$
with the measurable mappings $k_j\colon S_3 \to S_j$, $j=1,2$.

Now let $(X,\Xi)\defin(S_3,\Sigma_3)$
and for $x\in X$, $j=1,2$, define
\begin{align*}
  \mu^x_j &\defin \tau^{k_j(x)}_j  & j=1,2\\
  \mu^x_0 &\defin \tau^{h_1(k_1(x))}_0 = \tau^{h_2(k_2(x))}_0.
\end{align*}

Since the maps $k_j$, $h_j$ and $x\mapsto \tau^x_j$ are measurable, 
it follows that $\mu_j$ are subprobability kernels.
By Theorem~\ref{th:main}
there exists a semipullback $\mu_3$ in the category
of Radon subprobability kernels from $X=S_3$.
For $A\in\Sigma_j$, $j=1,2$, and $x\in S_3$ we have
\[
\mu^x_3(k^{-1}_j(A)) = \mu^x_j(A) =  \tau^{k_j(x)}_j(A),
\]
which means that $\mu_3$ is also a semipullback in the LMP category.
That concludes the proof for the case of a singleton label set $L$.

Now consider an arbitrary label set $L$. 
We have just proved that for each $a\in L$
there exists a semipullback $(S_3,\Sigma_3,\tau_a)$ in which 
$S_3$ and $\Sigma_3$ do not depend on $a$.
But that means that $(S_3,\Sigma_3,\{\tau_a : a\in L\})$
is a semipullback in the category of the LMP labelled by $L$.
\end{proof}

\subsection{Universally measurable labelled Markov processes}
\label{sec:univ-meas-label-mark-proc}

In Theorem~\ref{th:Markov} we assume that each measure $\tau_a(s,\cdot)$ is Radon.
It may be more convenient to have instead a single restriction on the underlying space $S$,
as in the next theorem.

\begin{defi}
  A measurable space  $(S,\Sigma)$ is a \emph{separable universally measurable
    space} if it is isomorphic to a  universally measurable subset of
  a separable completely metrizable (``Polish'') space with the trace
  of the Borel  $\sigma$-algebra.
\end{defi}

\begin{thm}
    \label{th:lmp-over-unvers-semipullb}
  Consider the category in which objects are LMP $(S,\Sigma, \{\tau_a : a\in L \})$
  such that $S$ is a separable universally measurable space,
  with zigzag morphisms.
  In this category, every cospan has a semipullback.
\end{thm}

\begin{proof}

Let
\[
  (S_1,\Sigma_1,\{\tau_{1a} : a\in L \})
  \rightarrow (S_0,\Sigma_0,\{\tau_{0a} : a\in L \})
  \leftarrow (S_2,\Sigma_2,\{\tau_{2a} : a\in L \})
\]
be a cospan of LMP with $S_j$ separable universally measurable spaces.
Then each $(S_j,\Sigma_j)$ is isomorphic to some
$(X_j,\B(Y_j)\restriction X_j)$ where $X_j$ is a universally measurable
subset of a Polish space $Y_j$ and $\B(Y_j)$ is its Borel
$\sigma$-algebra. Since $Y_j$ is a Radon space, by
\cite[434F(c)]{fremlinMTvol4} we conclude that every Borel measure
on $X_j$ is Radon. 

Let $(S_3,\Sigma_3,\{\tau_{3a} : a\in L \})$ be a semipullback 
with the properties from Theorem~\ref{th:Markov}.
In particular, $S_3$ is a measurable subset of $S_1\times S_2$.
It remains to prove that $S_3$ is a separable universally measurable space.
  There exists a  measurable isomorphism
  \[
  f:(S_1\times S_2,\Sigma_1\sigmagen\Sigma_2)\to
  (X_1\times X_2, \B(Y_1\times Y_2)\restriction
  X_1\times X_2),
  \]
  $X_1\times X_2$ is universally measurable by
  \cite[434X(c)]{fremlinMTvol4}. But then
  $(S_3,\Sigma_3)=(S_3,\Sigma_1\sigmagen\Sigma_2\restriction S_3)$ is
  isomorphic to $(f(S_3),\B(Y_1\times Y_2)\restriction f(S_3))$, where
  $f(S_3)\in\B(Y_1\times Y_2)\restriction X_1\times X_2$ since
  $S_3\in\Sigma_1\sigmagen\Sigma_2$.
\end{proof}

In \cite{Pedro20111048} it was asked whether behaviorally equivalent
LMP over coanalytic spaces were bisimilar. We can answer this question
affirmatively. Recall that a metric space is coanalytic if it
is homeomorphic to the complement of an analytic subset of a Polish space.
We say that a measurable space is
\emph{coanalytic} if it is isomorphic to the Borel space of a
coanalytic metric space.
\begin{cor}
  The category of LMP over coanalytic measurable spaces has
  semipullbacks.
\end{cor}
\begin{proof}
  This follows by essentially the same  argument for
  Theorem~\ref{th:lmp-over-unvers-semipullb}, showing that a cospan of
  coanalytic measurable spaces has a coanalytic pullback, and that
  coanalytic sets are universally measurable.
\end{proof}

In the same way, every analytic space with its Borel $\sigma$-algebra
is a separable universally measurable space; hence we obtain
Edalat's result \cite{Edalat} as a corollary to
Theorem~\ref{th:lmp-over-unvers-semipullb} as well.

\section{Counterexamples}
\label{sec:counterexamples}
The key assumption in previous sections is that each measure is
defined on the Borel $\sigma$-algebra. The results no longer hold without
that assumption, even for $\sigma$-algebras of subsets of $[0,1]$. The
counterexample in \cite{Pedro20111048} uses a $\sigma$-algebra larger
than the Borel 
$\sigma$-algebra on $[0,1]$; we hint at this construction below. In
the opposite direction, the following 
counterexample uses $\sigma$-algebras that are smaller but still large
enough to separate the points of $[0,1]$.

\begin{exa}\label{exm:1}
Consider $(S,\Sigma)$ to be the interval $[0,1]$ with the
countable-cocountable $\sigma$-algebra, and let $\mu_0:\Sigma\to
\{0,1\}$ be the probability measure such that 
\[
\mu_0(Q)=1 \iff Q\text{ has countable complement.}
\]
Take $V\defin [0,\tfrac{1}{2}]$. It is straightforward to check that for any
different $r_1,r_2\in (0,1)$, the following maps 
\[
\mu_i(Q) \defin
\begin{cases}
  \mu(Q) & Q\in\Sigma\\
  r_i &  Q \in \Sigma_V\setminus\Sigma \text{ and } Q\sm V \text{ is countable} \\
  1-r_i &  \text{otherwise,}
\end{cases}
\]
are probability measures that extend $\mu_0$ to
$\Sigma_V\defin\sigma(\Sigma\cup\{V\})$.

By using these probability spaces we can replicate the idea of
\cite[Thm.~12]{Pedro20111048} to obtain a cospan of LMP that can not be
completed to a commutative square. We now sketch the construction. Fix
$s_0\in S$ and define LMP $\mathbf{S}_i\defin  (S, \Sigma_i,\tau_i)$,
with $\Sigma_1 = \Sigma_2 \defin \Sigma_V$ and $\Sigma_0\defin \Sigma$,
and
\[
\tau_i(s,A) := \begin{cases} 1 &s\neq s_0
  \text{ and } s_0\in A \\
  \mu_i(A) & s = s_0 \\
  0 &\text{otherwise}
\end{cases}
\]
for $i=0,1,2$, every $s\in S$, and  $A$ in the corresponding $\sigma$-algebra.

The identity
maps $\mathit{Id}_S : \mathbf{S}_i \to \mathbf{S}_0$ form a cospan
of zigzags, and it can be seen that  there are no
$\mathbf{S}$ and zigzag maps $h_i: \mathbf{S}
\to \mathbf{S}_i$ ($i=1,2$) completing that cospan to a semipullback.
\end{exa}
\section*{Acknowledgment}
  \noindent
  We would like to thank Zolt\'an Vidny\'ansky for initial suggestions
  regarding 
  Example~\ref{exm:1}. The second author would also like to thank Prof.~S.M.~Srivastava
  and  Prof.~S.~Todor\v{c}evi\'c for enlightening discussions around this
  problem.

\bibliographystyle{alpha}
\bibliography{Semipullback_LMCS_typeset}

\end{document}